\documentclass{article}
\usepackage{amssymb,amsmath,amsthm,graphicx,cite}

\def\utr{\, \underline{\triangleright}\, }
\def\otr{\, \overline{\triangleright}\, }

\def\bar{\overline}

\newtheorem{theorem}{Theorem}

\newtheorem{proposition}[theorem]{Proposition}
\newtheorem{corollary}[theorem]{Corollary}

\theoremstyle{definition}
\newtheorem{example}{Example}

\newtheorem{definition}{Definition}
\newtheorem{remark}{Remark}

\title{\Large \textbf{Biquandle Power Brackets}}

\date{}

\author{Neslihan G\"ug\"umc\"u\footnote{Email:neslihangugumcu@iyte.edu.tr }
\and
Sam Nelson\footnote{Email: Sam.Nelson@cmc.edu. 
Partially supported by Simons Foundation collaboration grant 702597.}}

\begin{document}
\maketitle

\begin{abstract}
In this paper, we introduce biquandle power brackets, an infinite family
of invariants of oriented links containing the classical skein invariants
and the quandle and biquandle 2-cocycle invariants as special cases. Biquandle 
power brackets are generalizations of biquandle brackets in which the values 
of Kauffman states also depend on the biquandle colors they admit. We provide 
example computations and discuss the relationship between these new invariants 
and the previous cases.
\end{abstract}

\parbox{6in} {\textsc{Keywords:} Biquandle Invariants of Oriented Links, 
Biquandle Brackets, Biquandle Power Brackets

\smallskip

}

\section{\large\textbf{Introduction}}\label{I}

An intersection point of mathematical physics and low-dimensional topology,
\textit{quantum knot invariants} are a family of algebraic structures related
to oriented knots and links in $\mathbb{R}^3$ which have enjoyed great success
over the last few decades. With applications in theoretical physics, quantum 
computation, 3-manifold theory and beyond, the study of quantum knot 
invariants is vibrant and active research area. In this paper we introduce
a new infinite family of quantum knot invariants called \textit{biquandle
power brackets}.

A \textit{biquandle} is an algebraic structure consisting of a set with 
two operations satisfying axioms encoding the Reidemeister moves from knot 
theory analogously to the way the group axioms encode the properties of 
symmetry in geometry. Every 
oriented knot or link $L$ has an associated \textit{fundamental biquandle}
denoted $\mathcal{B}(L)$ whose isomorphism class is a strong invariant -- 
indeed, a complete invariant up to mirror image \cite{EN}. Every finite 
biquandle $X$ determines a multiset-valued invariant of oriented knots and 
links called the 
\textit{biquandle homset invariant}, denoted $\mathrm{Hom}(\mathcal{B}(L),X)$.
A choice of diagram for $L$ gives us a concrete representation for the
elements of the homset analogously to the way a choice of basis for a vector 
space gives us matrix representations of linear transformations; a different 
choice of diagram yields representations of the homset elements related to 
the original via $X$-colored Reidemeister moves, analogous to applying a 
change-of-basis matrix to a matrix encoding a linear transformation.

The cardinality of the homset is a natural number-valued invariant of 
oriented knots and links called  the \textit{biquandle counting invariant},
denoted $\Phi_X^{\mathbb{Z}}(X)$. An invariant $\phi$ of $X$-colored 
Reidemeister moves yields an \textit{enhancement} of the counting invariant,
a generally stronger invariant which specializes to $\Phi_X^{\mathbb{Z}}$. 
Introduced in \cite{NOR} and later studied in \cite{FN, GNO, HVW, IM, NO} etc., 
\textit{biquandle brackets} are quantum invariants of biquandle homset elements 
generalizing the Kauffman bracket/Jones polynomial to the case of 
biquandle-colored oriented knots and links. Biquandle brackets are defined 
via a state-sum with skein relations in which the 
skein coefficients are functions of the biquandle colors. In standard
biquandles brackets, the values of the components in the Kauffman states
are the same constant value $\delta=-A_{x,y}B_{x,y}^{-1}-A_{x,y}^{-1}B_{x,y}$. 
See \cite{NOR, NO} for more.

The set of biquandle bracket invariants is infinite, with infinite choices
for both coloring biquandle $X$ and coefficient ring $R$. Special cases include
the classical skein invariants such as the Alexander and Jones polynomials
as well as the quandle and biquandle 2-cocycle invariants. In \cite{FN}
biquandle brackets are categorified using biquandle coloring quivers, 
further strengthening this family of invariants.

In this paper we use the trace diagram formulation of biquandle brackets
from \cite{NO} to generalize the biquandle bracket definition, defining
\textit{biquandle power brackets} in which the
values of Kauffman state components also depend on the biquandle colors
they carry, with the previous biquandle brackets forming special cases 
with constant $\delta$ value. As an added bonus, we are able to consider
non-invertible skein coefficients, revealing some biquandle power brackets
with constant $\delta$ which do not fit the previous definition. The paper
is organized as follows. In Section \ref{BB} we review the basics of 
biquandles and the biquandle homset invariant. In Section \ref{BPB} we
introduce the new biquandle power brackets and prove that they define
invariants. In Section \ref{E} we collect some small examples to show that
the new invariants are proper enhancements, and we conclude with some 
questions for future research in Section \ref{Q}.

The second author would like to acknowledge the warm hospitality of IYTE
where the bulk of this paper was written during his sabbatical visit.

\section{\large\textbf{Biquandle Basics}}\label{BB}

We begin by reviewing preliminary notions of the theory of biquandles.

	\begin{definition}\label{def:biquandle}
		Let $X$ be a set that is endowed with two operations $\utr,\otr:X\times X\rightarrow X$. $(X,\utr,\otr)$ is called a \textit{biquandle} if for all $x,y,z\in X$, 
		\begin{enumerate}
		    
			\item $x\utr x=x\otr x$,
			\item the following maps are invertible,
			\begin{align*}
				&\alpha_y:X\rightarrow X, \qquad \beta_y:X\rightarrow X,\qquad S:X\times X\rightarrow X\times X\\
				&\quad x\mapsto x\otr y,\qquad x\mapsto x\utr y,\qquad\quad (x,y)\mapsto (y\otr x, x\utr y),
			\end{align*}
			\item the following exchange laws hold, \begin{align*}
				(x\utr y)\utr(z\utr y)=(x\utr z)\utr(y\otr z),\\
				(x\utr y)\otr(z\utr y)=(x\otr z)\utr(y\otr z),\\
				(x\otr y)\otr(z\otr y)=(x\otr z)\otr(y\utr z).
			\end{align*}
		\end{enumerate}
\end{definition}
\begin{example}
    \begin{enumerate}

 \item Let X be a set and $\sigma : X \rightarrow X$ be a bijection. We define $x\utr y = x \otr y =\sigma(x)$. Then, $(X, \utr, \otr)$ is called a \textit{constant action biquandle}. 
\item Let $R$ be a commutative ring with identity and $X$ be a module over $R[s^{\pm 1}, t^{\pm 1}]$. The \textit{Alexander biquandle} over $X$ is defined via
$x\utr y = tx + (s-t)y$ and $x \otr y = sx$.
\item For finite biquandle structures on the set $X=\{1,2,\dots, n\}$ we can
specify operation tables. For example, the smallest nontrivial biquandle
is $X=\{1,2\}$ with tables
\[
\begin{array}{r|rr}
\utr & 1 & 2 \\ \hline
1 & 2 & 2 \\ 
2 & 1 & 1
\end{array}\quad 
\begin{array}{r|rr}
\otr & 1 & 2 \\ \hline
1 & 2 & 2 \\ 
2 & 1 & 1
\end{array}.
\]
		
		\end{enumerate}
\end{example}

\begin{definition} 
Let $X, Y$ be biquandles. A \textit{biquandle homomorphism} is a map 
$f:X \rightarrow Y$ such that for all $x, y\in X$, 
$f(x\utr y) = f(x)\utr f(y)$ and $f(x\otr y) = f(x)\otr f(y)$.
\end{definition}


We can color an oriented link diagram with elements of a biquandle. 
A \textit{semiarc} of a link diagram is a portion of the diagram that lies between two adjacent crossings. To color the diagram with a biquandle $X$, an element of $X$ is assigned to each semiarc such that at each crossing, the relations
\[\includegraphics{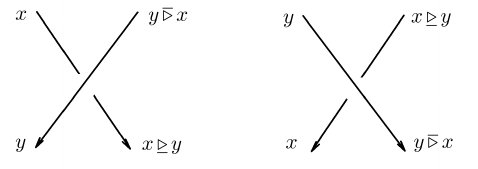}\]
are satisfied. The axioms of a biquandle are induced by oriented Reidemeister moves. Precisely, the axiom set of a biquandle consists of the conditions that provide invariance on the biquandle coloring of an oriented link diagram under oriented Reidemeister moves. We verify below that the exchange laws (Axiom 3)  given in Definition \ref{def:biquandle} are induced by the following oriented Reidemeister three move that take place on a biquandle colored link diagram. The verification for the remaining axioms are left to the reader.

\[\includegraphics{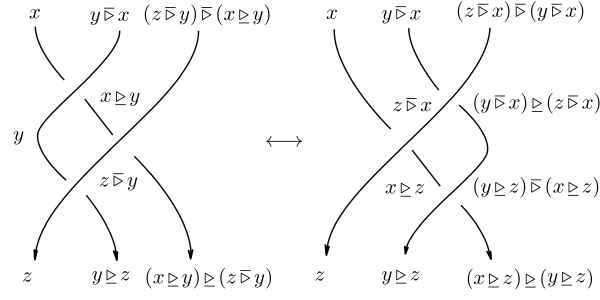}\]

\begin{definition}
Let $L$ be an oriented link diagram whose semi-arcs are labeled. 
Then the fundamental biquandle of $L$, denoted $\mathcal{B}(L)$,  
is the set of all equivalence classes of finite words
generated by the labels of semi-arcs of $L$ and the two operations (together 
with symbols required by the inverse axiom such as $\alpha^{-1}(x_j)$ etc.)
modulo the relations on the labels induced at each crossing 
\[\includegraphics{ng-sn2-15.pdf}\]
plus the biquandle axioms.
\end{definition}

It is easy to check from the axioms (see \cite{EN} for example) that:
\begin{theorem}
The fundamental biquandle is an invariant of oriented links.
\end{theorem}

\begin{remark}
In fact, more can be said depending on the category of knotted object in 
question. For oriented knots, the fundamental biquandle determines the 
knot quandle, which (by Joyce \cite{J}) is a complete invariant up to ambient
homeomorphism. For virtual knots, the fundamental biquandle is conjectured to
be a complete invariant up to orientation-reversed mirror image, but this is
still open.
\end{remark}

Each coloring $c$ of an oriented link diagram $L$ with a finite biquandle 
$X$ corresponds to a unique well-defined map $c : \mathcal{B}(L) \rightarrow X$ 
defined as $c(x_i) = X_j$ where $x_i$ is a generator of $\mathcal{B}(L)$ 
and $X_j$ is the element of the biquandle $X$ coloring the semiarc whose 
generator is $x_i$. It is clear that $c(x_i\utr x_j) = c(x_i)\utr c(x_j)$ 
and $c(x_i\otr x_j) = c(x_i)\otr c(x_j)$ and so each coloring corresponds 
to a unique homomorphism. The set of all homomorphisms that is induced by a 
biquandle coloring of $L$ with $X$  is denoted by 
$\mathrm{Hom}(\mathcal{B}(L), X)$. It then follows (see also \cite{EN}) that:

\begin{theorem} 
Let $L$ denote an oriented link diagram and $X$ be a finite biquandle.
The set of homomorphisms $\mathrm{Hom}(\mathcal{B}(L), X)$ is an invariant 
of oriented links. 
\end{theorem}

\begin{definition}
The cardinality of the biquandle homset, deoted $\Phi_X^{\mathbb{Z}}(L)$, 
is a nonnegative integer-valued invariant of oreinted links called the 
\textit{biquandle counting invariant}.
\end{definition}

\begin{remark}
A biquandle coloring of an oriented link diagram by a finite biquandle $X$
represents an element of the biquandle homset $\mathrm{Hom}(\mathcal{B}(L),X)$ 
analogously to a matrix representing a linear transformation between vector 
spaces. Changing an $X$-colored diagram by $X$-colored Reidemeister moves 
changes the representation of the biquandle homomorphism analogously to 
conjugating a matrix by a change-of-basis matrix.
\end{remark}

\section{\large\textbf{Biquandle Power Brackets}}\label{BPB}

We begin with our main definition.

\begin{definition}\label{def:1}
Let $X$ be a finite biquandle and $R$ a commutative ring with identity.
A \textit{biquandle power bracket} on $X$ with coefficients in $R$ consists
of an invertible constant $w\in R$, 
a map $\delta:\mathcal{P}(X)\to R$ from the power set of $X$ to $R$ and
a quadruple of maps $A,B,\bar{A},\bar{B}:X\times X\to R$ from ordered pairs 
of elements of $X$ to $R$ satisfying the following axioms:
\begin{itemize}
\item[(i)] For all $C\in \mathcal{P}(X)$ and all $x\in C$, 
\begin{eqnarray*}
w\delta(C) 
& = & A_{x,x}\delta(C)\delta(\{x\utr x\})+B_{x,x}\delta(C\cup\{x\utr x\}) \\
w^{-1} \delta(C) 
& = & \bar{A}_{x,x}\delta(C)\delta(\{x\utr x\})
+\bar{B}_{x,x}\delta(C\cup\{x\utr x\}),
\end{eqnarray*}
and for all $x$ such that $x\utr x\in C$,
\begin{eqnarray*}
w\delta(C)  
& = & A_{x,x}\delta(C)\delta(\{x\})+B_{x,x}\delta(C\cup\{x\}) \\
w^{-1} \delta(C) 
& = & \overline{A}_{x,x}\delta(C)\delta(\{x\})
+\overline{B}_{x,x}\delta(C\cup\{x\}), 
\end{eqnarray*}
\item[(ii.i)]
For all $C_1,C_2\in \mathcal{P}(X)$, if $x,y\in C_1\cap C_2$ then
\begin{eqnarray*}
\delta(C_1\cup C_2)
& = & A_{x,y}\bar{A}_{x,y}\delta(C_1)\delta(C_2)\delta(\{x\utr y,y\otr x\})
+B_{x,y}\bar{A}_{x,y}\delta(C_1\cup\{x\utr y,y\otr x\})\delta(C_2)\\ & & 
+A_{x,y}\bar{B}_{x,y}\delta(C_1)\delta(C_2\cup\{x\utr y,y\otr x\})  
+B_{x,y}\bar{B}_{x,y}\delta(C_1\cup C_2\cup\{x\utr y,y\otr x\}),
\end{eqnarray*} 
and if $x\in C_1$  and $y\in C_2$ then
\begin{eqnarray*}
\delta(C_1)\delta(C_2)
& = & A_{x,y}\bar{A}_{x,y}\delta(C_1\cup C_2)\delta(\{x\utr y,y\otr x\})
+B_{x,y}\bar{A}_{x,y}\delta(C_1\cup C_2\cup\{x\utr y,y\otr x\})\\ & & 
+A_{x,y}\bar{B}_{x,y}\delta(C_1\cup C_2\cup\{x\utr y,y\otr x\}) 
+B_{x,y}\bar{B}_{x,y}\delta(C_1\cup\{y\otr x\})\delta(C_2\cup\{x\utr y\}),
\end{eqnarray*}

\item[(ii.ii)]
For all $C_1,C_2\in \mathcal{P}(X)$, 
for all $x,y\in X$ such that $y\otr x,x\utr y\in C_1\cup C_2$, we have
\begin{eqnarray*}
\delta(C_1\cup C_2) 
& = & A_{x,y}\bar{A}_{x,y}\delta(C_1)\delta(C_2)\delta(\{x,y\})
+A_{x,y}\bar{B}_{x,y}\delta(C_1\cup \{x,y\})\delta(C_2) \\ & &
+B_{x,y}\bar{A}_{x,y}\delta(C_1)\delta(C_2\cup\{x,y\})
+B_{x,y}\bar{B}_{x,y}\delta(C_1\cup C_2\cup\{x,y\}),
\end{eqnarray*}
and for all $x,y\in X$ such that $y\otr x\in C_1$ and $x\utr y\in C_2$ we have
\begin{eqnarray*}
\delta(C_1)\delta(C_2) 
& = & A_{x,y}\bar{A}_{x,y}\delta(C_1\cup C_2)\delta(\{x,y\})
+A_{x,y}\bar{B}_{x,y}\delta(C_1\cup C_2 \cup \{x,y\}) \\ & &
+B_{x,y}\bar{A}_{x,y}\delta(C_1\cup C_2 \cup \{x,y\})
+B_{x,y}\bar{B}_{x,y}\delta(C_1\cup \{x\})\delta(C_2\cup\{y\})
\end{eqnarray*}
\item[(iii.i)]
For all $C_1,C_2,C_3\in \mathcal{P}(X)$ and all $x,y,z\in X$ such that
$x,z\in C_1$, $y\otr x,y\utr z\in C_2$ and 
$(z\otr x)\otr (y\otr x), (x\utr z)\utr(y\utr z)\in C_3$ we have

\hskip -0.5 in
\scalebox{0.85}{$\begin{array}{rcl}
\begin{array}{r}
A_{x,y}A_{y,z}A_{x\utr y,z\otr y}\delta(C_1\cup \{y\})\delta(C_2\cup\{x\utr y,z\otr y\})\delta(C_3)  \\
+A_{x,y}A_{y,z}B_{x\utr y,z\otr y}\delta(C_1\cup \{y\})\delta(C_2\cup C_3\cup \{x\utr y,z\otr y\})  \\
+A_{x,y}B_{y,z}A_{x\utr y,z\otr y}\delta(C_1\cup C_2\cup \{y, x\utr y,z\otr y\})\delta(C_3)  \\
+A_{x,y}B_{y,z}B_{x\utr y,z\otr y}\delta(C_1\cup C_2\cup C_3\cup \{y, x\utr y,z\otr y\})  \\
+B_{x,y}A_{y,z}A_{x\utr y,z\otr y}\delta(C_1\cup C_2\cup \{y, x\utr y,z\otr y\})\delta(C_3)  \\
+B_{x,y}A_{y,z}B_{x\utr y,z\otr y}\delta(C_1\cup C_2\cup C_3\cup \{y, x\utr y,z\otr y\})  \\
+B_{x,y}B_{y,z}A_{x\utr y,z\otr y}\delta(C_1\cup C_2)\delta(C_3)\delta(\{y, x\utr y,z\otr y\})  \\
+B_{x,y}B_{y,z}B_{x\utr y,z\otr y}\delta(C_1\cup C_2)\delta(C_3 \cup \{y, x\utr y,z\otr y\}) 
\end{array}
& = & 
\begin{array}{l}
A_{x,z}A_{y\otr x, z\otr x}A_{x\utr z,y\utr z}\delta(C_1)\delta(C_2\cup\{z\otr x,x\utr z\})\delta(C_3\cup\{(y\utr z)\otr(x\utr z)\}) \\
+A_{x,z}A_{y\otr x, z\otr x}B_{x\utr z,y\utr z}\delta(C_1)\delta(C_2\cup C_3\cup\{z\otr x,x\utr z,(y\utr z)\otr(x\utr z)\}) \\
+A_{x,z}B_{y\otr x, z\otr x}A_{x\utr z,y\utr z}\delta(C_1)\delta(C_2\cup C_3\cup\{z\otr x,x\utr z,(y\utr z)\otr(x\utr z)\}) \\
+A_{x,z}B_{y\otr x, z\otr x}B_{x\utr z,y\utr z}\delta(C_1)\delta(C_2\cup C_3)\delta(\{z\otr x,x\utr z,(y\utr z)\otr(x\utr z)\}) \\
+B_{x,z}A_{y\otr x, z\otr x}A_{x\utr z,y\utr z}\delta(C_1\cup C_2\cup\{z\otr x,x\utr z\})\delta(C_3\cup\{(y\utr z)\otr(x\utr z)\}) \\
+B_{x,z}A_{y\otr x, z\otr x}B_{x\utr z,y\utr z}\delta(C_1\cup C_2\cup C_3\cup\{z\otr x,x\utr z,(y\utr z)\otr(x\utr z)\}) \\
+B_{x,z}B_{y\otr x, z\otr x}A_{x\utr z,y\utr z}\delta(C_1\cup C_2\cup C_3\cup\{z\otr x,x\utr z,(y\utr z)\otr(x\utr z)\}) \\
+B_{x,z}B_{y\otr x, z\otr x}B_{x\utr z,y\utr z}\delta(C_1\cup\{z\otr x,x\utr z,(y\utr z)\otr(x\utr z)\})\delta(C_2\cup C_3) 
\end{array}
\end{array}$}
\item[(iii.ii)]
For all $C_1,C_2,C_3\in \mathcal{P}(X)$ and all $x,y,z\in X$ such that
$x,y\otr x\in C_1$, $(z\otr x)\otr (y\otr x), (x\utr z)\utr(y\utr z)\in C_2$ 
and $z,y\utr z\in C_3$ we have

\hskip -0.5 in
\scalebox{0.85}{$\begin{array}{rcl}
\begin{array}{r}
A_{x,y}A_{y,z}A_{x\utr y,z\otr y}\delta(C_1\cup C_3\cup \{y, x\utr y,z\otr y\})\delta(C_2)  \\
+A_{x,y}A_{y,z}B_{x\utr y,z\otr y}\delta(C_1\cup C_2\cup C_3\cup \{y,x\utr y,z\otr y\})  \\
+A_{x,y}B_{y,z}A_{x\utr y,z\otr y}\delta(C_1\cup \{y, x\utr y,z\otr y\})\delta(C_2)\delta(C_3)  \\
+A_{x,y}B_{y,z}B_{x\utr y,z\otr y}\delta(C_1\cup C_2\cup \{y, x\utr y,z\otr y\})\delta(C_3)  \\
+B_{x,y}A_{y,z}A_{x\utr y,z\otr y}\delta(C_1)\delta(C_2)\delta(C_3\cup \{y, x\utr y,z\otr y\})  \\
+B_{x,y}A_{y,z}B_{x\utr y,z\otr y}\delta(C_1)\delta(C_2\cup C_3\cup \{y, x\utr y,z\otr y\})  \\
+B_{x,y}B_{y,z}A_{x\utr y,z\otr y}\delta(C_1)\delta (C_2)\delta(C_3)\delta(\{y, x\utr y,z\otr y\})  \\
+B_{x,y}B_{y,z}B_{x\utr y,z\otr y}\delta(C_1)\delta(C_2\cup \{y, x\utr y,z\otr y\})\delta(C_3) 
\end{array}
& = & 
\begin{array}{l}
A_{x,z}A_{y\otr x, z\otr x}A_{x\utr z,y\utr z}\delta(C_1\cup C_3\cup\{z\otr x,x\utr z\})\delta(C_2\cup\{(y\utr z)\otr(x\utr z)\}) \\
+A_{x,z}A_{y\otr x, z\otr x}B_{x\utr z,y\utr z}\delta(C_1\cup C_2\cup C_3\cup\{z\otr x,x\utr z,(y\utr z)\otr(x\utr z)\}) \\
+A_{x,z}B_{y\otr x, z\otr x}A_{x\utr z,y\utr z}\delta(C_1\cup C_2\cup C_3\cup\{z\otr x,x\utr z,(y\utr z)\otr(x\utr z)\}) \\
+A_{x,z}B_{y\otr x, z\otr x}B_{x\utr z,y\utr z}\delta(C_1\cup C_2\cup C_3)\delta(\{z\otr x,x\utr z,(y\utr z)\otr(x\utr z)\}) \\
+B_{x,z}A_{y\otr x, z\otr x}A_{x\utr z,y\utr z}\delta(C_1\cup \{z\otr x\})\delta(C_2\cup \{(y\utr z)\otr(x\utr z)\})\delta(C_3\cup\{x\utr z\}) \\
+B_{x,z}A_{y\otr x, z\otr x}B_{x\utr z,y\utr z}\delta(C_1\cup\{z\otr x\})\delta(C_2\cup C_3\cup\{x\utr z,(y\utr z)\otr(x\utr z)\}) \\
+B_{x,z}B_{y\otr x, z\otr x}A_{x\utr z,y\utr z}\delta(C_1\cup C_2\cup\{z\otr x,(y\utr z)\otr(x\utr z)\})\delta(C_3\cup\{x\utr z\}) \\
+B_{x,z}B_{y\otr x, z\otr x}B_{x\utr z,y\utr z}\delta(C_1\cup C_2\cup C_3 \cup\{z\otr x,x\utr z,(y\utr z)\otr(x\utr z)\})
\end{array}
\end{array}$}
\item[(iii.iii)]
For all $C_1,C_2,C_3\in \mathcal{P}(X)$ and all $x,y,z\in X$ such that
$x,z\in C_1$, $y\otr x,(z\otr x)\otr (y\otr x)\in C_2$ and 
$y\utr z, (x\utr z)\utr(y\utr z)\in C_3$ we have

\hskip -0.5 in
\scalebox{0.85}{$\begin{array}{rcl}
\begin{array}{r}
A_{x,y}A_{y,z}A_{x\utr y,z\otr y}\delta(C_1\cup \{y\})\delta(C_2\cup C_3\cup \{x\utr y,z\otr y\}) \\
+A_{x,y}A_{y,z}B_{x\utr y,z\otr y}\delta(C_1\cup \{y\})\delta(C_2\cup \{x\utr y\})\delta(C_3\cup\{z\otr y\}) \\
+A_{x,y}B_{y,z}A_{x\utr y,z\otr y}\delta(C_1\cup C_2\cup C_3\cup\{y, x\utr y,z\otr y\})  \\
+A_{x,y}B_{y,z}B_{x\utr y,z\otr y}\delta(C_1\cup C_3\cup\{y, z\otr y\})\delta(C_2\cup\{x\utr y\})  \\
+B_{x,y}A_{y,z}A_{x\utr y,z\otr y}\delta(C_1\cup C_2\cup C_3\cup\{y, x\utr y,z\otr y\})  \\
+B_{x,y}A_{y,z}B_{x\utr y,z\otr y}\delta(C_1\cup C_2\cup \{y,x\utr y\})\delta(C_3\cup \{z\otr y\})  \\
+B_{x,y}B_{y,z}A_{x\utr y,z\otr y}\delta(C_1\cup C_2\cup C_3)\delta(\{y, x\utr y,z\otr y\})  \\
+B_{x,y}B_{y,z}B_{x\utr y,z\otr y}\delta(C_1\cup C_2\cup C_3\cup \{y, x\utr y,z\otr y\}) 
\end{array}
& = & 
\begin{array}{l}
A_{x,z}A_{y\otr x, z\otr x}A_{x\utr z,y\utr z}\delta(C_1)\delta(C_2 \cup C_3\cup\{z\otr x,x\utr z, (y\utr z)\otr(x\utr z)\}) \\
+A_{x,z}A_{y\otr x, z\otr x}B_{x\utr z,y\utr z}\delta(C_1)\delta(C_2 \cup\{z\otr x,x\utr z, (y\utr z)\otr(x\utr z)\})\delta(C_3) \\
+A_{x,z}B_{y\otr x, z\otr x}A_{x\utr z,y\utr z}\delta(C_1)\delta(C_2) \delta(C_3\cup\{z\otr x,x\utr z, (y\utr z)\otr(x\utr z)\})\\
+A_{x,z}B_{y\otr x, z\otr x}B_{x\utr z,y\utr z}\delta(C_1)\delta(C_2) \delta(C_3)\delta(\{z\otr x,x\utr z, (y\utr z)\otr(x\utr z)\}) \\
+B_{x,z}A_{y\otr x, z\otr x}A_{x\utr z,y\utr z}\delta(C_1\cup C_2\cup C_3\cup \{z\otr x,x\utr z\,(y\utr z)\otr(x\utr z)\}) \\
+B_{x,z}A_{y\otr x, z\otr x}B_{x\utr z,y\utr z}\delta(C_1\cup C_2\cup \{z\otr x,x\utr z,(y\utr z)\otr(x\utr z)\})\delta(C_3) \\
+B_{x,z}B_{y\otr x, z\otr x}A_{x\utr z,y\utr z}\delta(C_1\cup C_3\cup \{z\otr x,x\utr z,(y\utr z)\otr(x\utr z)\})\delta(C_2) \\
+B_{x,z}B_{y\otr x, z\otr x}B_{x\utr z,y\utr z}\delta(C_1\cup\{z\otr x,x\utr z,(y\utr z)\otr(x\utr z)\})\delta(C_2)\delta(C_3)
\end{array}
\end{array}$}
\item[(iii.iv)]
For all $C_1,C_2,C_3\in \mathcal{P}(X)$ and all $x,y,z\in X$ such that
$x,(x\utr z)\utr(y\utr z)\in C_1$, $y\otr x,(z\otr x)\otr (y\otr x)\in C_2$ and 
$z,y\utr z\in C_3$ we have

\hskip -0.5 in
\scalebox{0.85}{$\begin{array}{rcl}
\begin{array}{r}
A_{x,y}A_{y,z}A_{x\utr y,z\otr y}\delta(C_1\cup C_2\cup C_3\cup \{y,x\utr y,z\otr y\}) \\
+A_{x,y}A_{y,z}B_{x\utr y,z\otr y}\delta(C_1\cup C_3\cup \{y,z\otr y\})\delta(C_2\cup \{x\utr y\}) \\
+A_{x,y}B_{y,z}A_{x\utr y,z\otr y}\delta(C_1\cup C_2\cup\{y, x\utr y,z\otr y\})\delta(C_3)  \\
+A_{x,y}B_{y,z}B_{x\utr y,z\otr y}\delta(C_1\cup\{y,z\otr y\})\delta(C_2\cup\{x\utr y\})\delta(C_3)  \\
+B_{x,y}A_{y,z}A_{x\utr y,z\otr y}\delta(C_1\cup C_2)\delta(C_3\cup\{y, x\utr y,z\otr y\})  \\
+B_{x,y}A_{y,z}B_{x\utr y,z\otr y}\delta(C_1\cup C_2\cup C_3\cup\{y, x\utr y,z\otr y\})  \\
+B_{x,y}B_{y,z}A_{x\utr y,z\otr y}\delta(C_1\cup C_2)\delta(C_3)\delta(\{y, x\utr y,z\otr y\})  \\
+B_{x,y}B_{y,z}B_{x\utr y,z\otr y}\delta(C_1\cup C_2\cup \{y, x\utr y,z\otr y\})\delta(C_3) 
\end{array}
& = & 
\begin{array}{l}
A_{x,z}A_{y\otr x, z\otr x}A_{x\utr z,y\utr z}\delta(C_1\cup C_2 \cup C_3\cup\{z\otr x, x\utr z, (y\utr z)\otr(x\utr z)\}) \\
+A_{x,z}A_{y\otr x, z\otr x}B_{x\utr z,y\utr z}\delta(C_1\cup C_3)\delta(C_2\cup\{z\otr x, x\utr z, (y\utr z)\otr(x\utr z)\}) \\
+A_{x,z}B_{y\otr x, z\otr x}A_{x\utr z,y\utr z}\delta(C_1 \cup C_3\cup\{z\otr x,x\utr z, (y\utr z)\otr(x\utr z)\})\delta(C_2) \\
+A_{x,z}B_{y\otr x, z\otr x}B_{x\utr z,y\utr z}\delta(C_1\cup C_3)\delta(C_2)\delta(\{z\otr x, x\utr z, (y\utr z)\otr(x\utr z)\}) \\
+B_{x,z}A_{y\otr x, z\otr x}A_{x\utr z,y\utr z}\delta(C_1\cup C_2\cup \{z\otr x,(y\utr z)\otr(x\utr z)\})\delta(C_3\cup\{x\utr z\}) \\
+B_{x,z}A_{y\otr x, z\otr x}B_{x\utr z,y\utr z}\delta(C_1\cup C_2\cup C_3\cup \{z\otr x, x\utr z,(y\utr z)\otr(x\utr z)\}) \\
+B_{x,z}B_{y\otr x, z\otr x}A_{x\utr z,y\utr z}\delta(C_1\cup \{z\otr x,(y\utr z)\otr(x\utr z)\})\delta(C_2)\delta(C_3\cup\{x\utr z\}) \\
+B_{x,z}B_{y\otr x, z\otr x}B_{x\utr z,y\utr z}\delta(C_1\cup C_3\cup \{z\otr x, x\utr z,(y\utr z)\otr(x\utr z)\})\delta(C_2)
\end{array}
\end{array}$}
\item[(iii.v)]
For all $C_1,C_2,C_3\in \mathcal{P}(X)$ and all $x,y,z\in X$ such that
$x,y\otr x\in C_1$, $z,(z\otr x)\otr (y\otr x)\in C_2$ and 
$y\utr z, (x\utr z)\utr(y\utr z) \in C_3$ we have

\hskip -0.5 in
\scalebox{0.85}{$\begin{array}{rcl}
\begin{array}{r}
A_{x,y}A_{y,z}A_{x\utr y,z\otr y}\delta(C_1\cup C_2\cup C_3\cup \{y,x\utr y,z\otr y\}) \\
+A_{x,y}A_{y,z}B_{x\utr y,z\otr y}\delta(C_1\cup C_2\cup \{y,x\utr y\})\delta(C_3\cup \{z\otr y\}) \\
+A_{x,y}B_{y,z}A_{x\utr y,z\otr y}\delta(C_1\cup\{y, x\utr y,z\otr y\})\delta(C_2\cup C_3)  \\
+A_{x,y}B_{y,z}B_{x\utr y,z\otr y}\delta(C_1\cup C_2\cup C_3\cup \{y, x\utr y,z\otr y\})  \\
+B_{x,y}A_{y,z}A_{x\utr y,z\otr y}\delta(C_1)\delta(C_2\cup C_3\cup\{y, x\utr y,z\otr y\})  \\
+B_{x,y}A_{y,z}B_{x\utr y,z\otr y}\delta(C_1)\delta(C_2\cup \{y, x\utr y\})\delta(C_3\cup\{z\otr y\})  \\
+B_{x,y}B_{y,z}A_{x\utr y,z\otr y}\delta(C_1)\delta(C_2\cup C_3)\delta(\{y, x\utr y,z\otr y\})  \\
+B_{x,y}B_{y,z}B_{x\utr y,z\otr y}\delta(C_1)\delta(C_2\cup C_3\cup\{y, x\utr y,z\otr y\})  \\
\end{array}
& = & 
\begin{array}{l}
A_{x,z}A_{y\otr x, z\otr x}A_{x\utr z,y\utr z}\delta(C_1\cup C_2 \cup C_3\cup\{z\otr x,x\utr z, (y\utr z)\otr(x\utr z)\}) \\
+A_{x,z}A_{y\otr x, z\otr x}B_{x\utr z,y\utr z}\delta(C_1\cup C_2\cup\{z\otr x,x\utr z, (y\utr z)\otr(x\utr z)\})\delta(C_3) \\
+A_{x,z}B_{y\otr x, z\otr x}A_{x\utr z,y\utr z}\delta(C_1\cup C_2)\delta(C_3\cup\{z\otr x,x\utr z, (y\utr z)\otr(x\utr z)\}) \\
+A_{x,z}B_{y\otr x, z\otr x}B_{x\utr z,y\utr z}\delta(C_1\cup C_2)\delta(C_3)\delta(\{z\otr x,x\utr z, (y\utr z)\otr(x\utr z)\}) \\
+B_{x,z}A_{y\otr x, z\otr x}A_{x\utr z,y\utr z}\delta(C_1\cup \{z\otr x\})\delta(C_2\cup C_3\cup \{x\utr z, (y\utr z)\otr(x\utr z)\}) \\
+B_{x,z}A_{y\otr x, z\otr x}B_{x\utr z,y\utr z}\delta(C_1\cup \{z\otr x\})\delta(C_2\cup \{x\utr z,(y\utr z)\otr(x\utr z)\})\delta(C_3) \\
+B_{x,z}B_{y\otr x, z\otr x}A_{x\utr z,y\utr z}\delta(C_1\cup C_2\cup C_3\cup \{z\otr x,x\utr z,(y\utr z)\otr(x\utr z)\}) \\
+B_{x,z}B_{y\otr x, z\otr x}B_{x\utr z,y\utr z}\delta(C_1\cup C_2\cup \{z\otr x,x\utr z,(y\utr z)\otr(x\utr z)\})\delta(C_3)
\end{array}
\end{array}$}
\end{itemize}
\end{definition}

This definition is motivated by the following system of skein relations,
Kauffman state values and writhe-correction factors
indexed by sets of biquandle elements, generalizing biquandle brackets as 
defined originally in \cite{NOR} using the trace diagram formulation from
\cite{NO}.
\[\scalebox{0.9}{\includegraphics{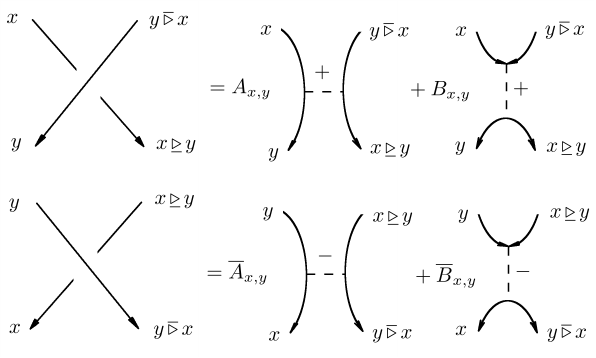}}\]
\[\scalebox{0.9}{\includegraphics{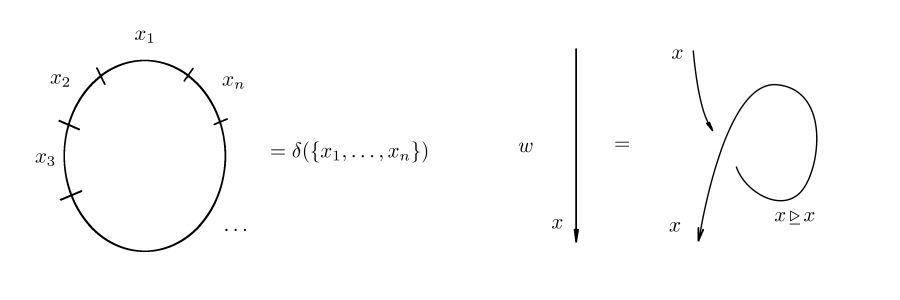}}\]

Indeed, the axioms in Definition \ref{def:1} are precisely the conditions 
required for the writhe-corrected state-sum of products of skein 
coefficients times state values, defined as products of delta-values of 
sets of biquandle colors on 
each component, to be invariant under the generating set of $X$-colored 
oriented Reidemeister moves including all four type I moves, both reverse 
type II moves and the all-positive crossing type III move; this set is 
shown to be a generating set of oriented Reidemeister moves in \cite{P}.

In the state-sum computation, we sum over all combinations of smoothings the 
products of skein coefficients times products of delta-values for the component
curves in the states. For a given Reidemeister move involving crossings
$c_1,\dots, c_k$, let us fix a choice of smoothing for each of the crossings 
$c_{k+1},\dots, c_n$ not involved in the move. The product of skein coefficients 
for $c_{k+1},\dots,c_n$ times delta values for the smoothed components not 
containing $c_1,\dots, c_k$ is then constant over the set of smoothings for
$c_1,\dots, c_k$ and thus factors out of the contributions to the overall 
state-sum for these states; it follows that for invariance we need equality 
of the sums of skein coefficients for $c_1,\dots, c_k$ times the product of 
delta values for the component curves for each crossingless way of closing 
the tangle containing $c_1,\dots, c_k$ on the two sides of the move.

For Reidemeister I moves, we have
\[\scalebox{0.9}{\includegraphics{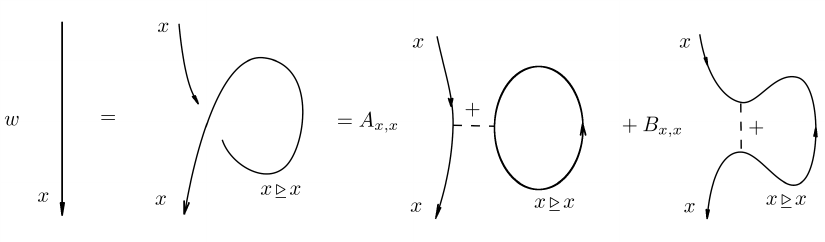}}\]
letting $C$ be the set of colors on the smoothed state component containing 
$x$ before introducing the crossing; then after the move we must have
\[w\delta(C)=A_{x,x}\delta(C)\delta\{x\utr x\}+B_{x,x}\delta(C\cup\{x\utr x\}),\] 
and similarly we have
\[\scalebox{0.9}{\includegraphics{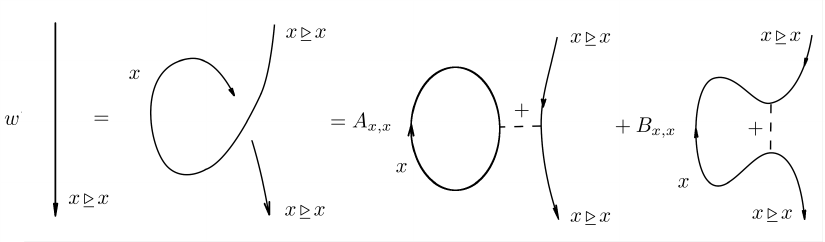}}\]
yielding 
\[w\delta(C)=A_{x,x}\delta(C)\delta(\{x\})+B_{x,x}\delta(C\cup\{x\})\] 
when $x\utr x\in C$. The negative-crossing counterparts yield the same 
equations with bars over the skein coefficients and $w$ replaced with $w^{-1}$.

For Reidemeister II moves, there are several cases. More precisely, 
in our generating set there are 
two oriented Reidemeister II moves, and for each we need to consider the 
possible connectivity of the components in the trace diagram Kauffman states. 
 
The reverse type II move with single-component closure $C=C_1\cup C_2$
in the Kauffman state with $x,y\in C_1\cup C_2$ 
\[\scalebox{0.8}{\includegraphics{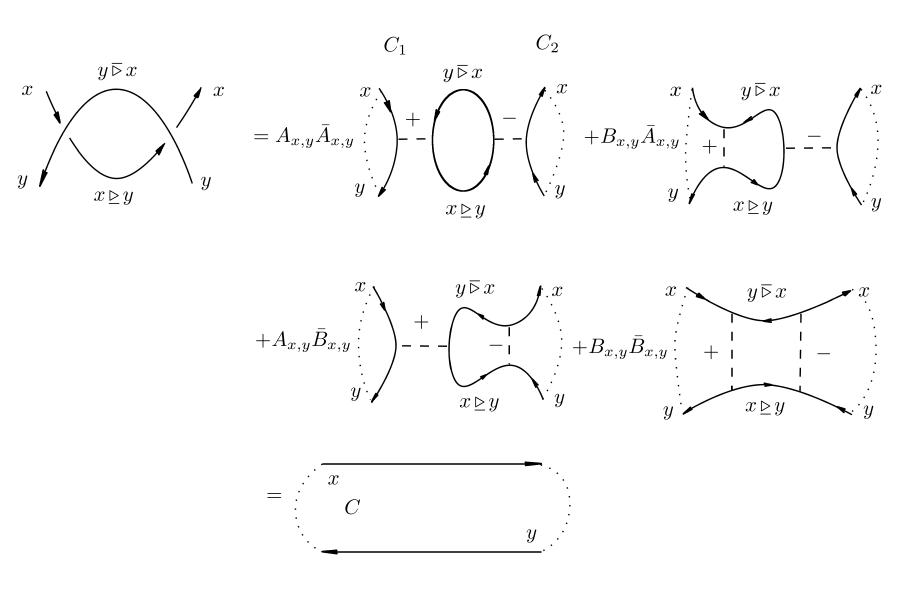}}\]
yields the requirement that
\begin{eqnarray*}
\delta(C_1\cup C_2)
& = & A_{x,y}\bar{A}_{x,y}\delta(C_1)\delta(C_2)\delta(\{x\utr y,y\otr x\})
+B_{x,y}\bar{A}_{x,y}\delta(C_1\cup\{x\utr y,y\otr x\})\delta(C_2)\\ & & 
+A_{x,y}\bar{B}_{x,y}\delta(C_1)\delta(C_2\cup\{x\utr y,y\otr x\})  
+B_{x,y}\bar{B}_{x,y}\delta(C_1\cup C_2\cup\{x\utr y,y\otr x\}),
\end{eqnarray*}
while the two-component closure with $x\in C_1$ and $y\in C_2$ 
\[\scalebox{0.75}{\includegraphics{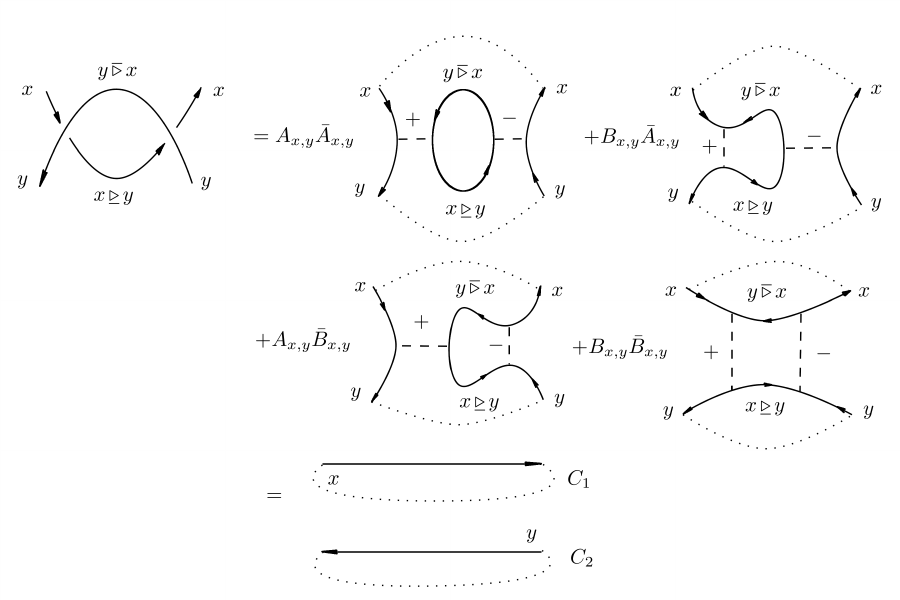}}\]
yields the requirement that 
\begin{eqnarray*}
\delta(C_1)\delta(C_2)
& = & A_{x,y}\bar{A}_{x,y}\delta(C_1\cup C_2)\delta(\{x\utr y,y\otr x\})
+B_{x,y}\bar{A}_{x,y}\delta(C_1\cup C_2\cup\{x\utr y,y\otr x\})\\ & & 
+A_{x,y}\bar{B}_{x,y}\delta(C_1\cup C_2\cup\{x\utr y,y\otr x\}) 
+B_{x,y}\bar{B}_{x,y}\delta(C_1\cup\{y\otr x\})\delta(C_2\cup\{x\utr y\}).
\end{eqnarray*}

The reverse type II move with single-component closure in the Kauffman 
state $C=C_1\cup C_2$ with $x\utr ,y\otr x\in C_1\cap C_2$
\[\scalebox{0.8}{\includegraphics{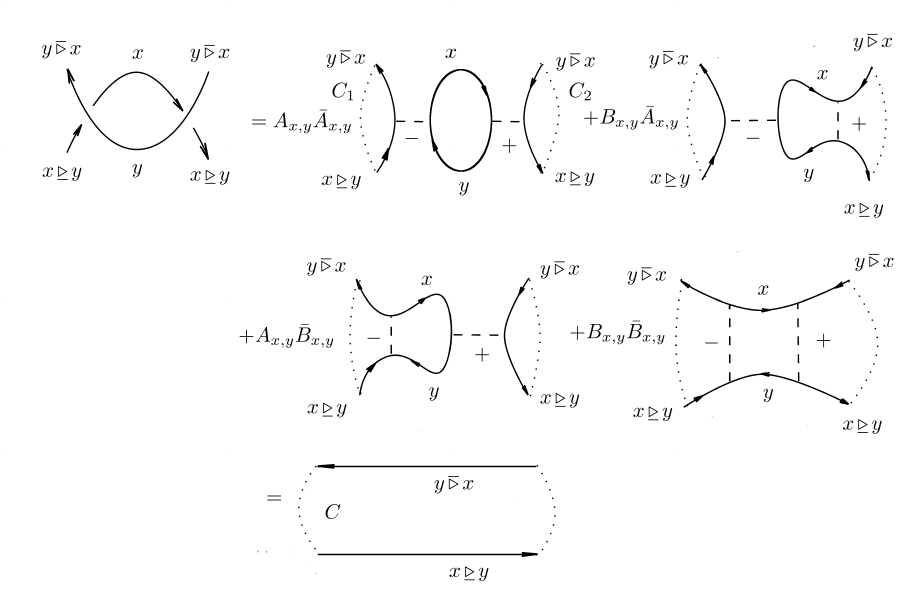}}\]
yields the requirement that 
\begin{eqnarray*}
\delta(C_1\cup C_2)
& = & A_{x,y}\bar{A}_{x,y}\delta(C_1)\delta(C_2)\delta(\{x,y\})
+B_{x,y}\bar{A}_{x,y}\delta(C_1)\delta(C_2\cup\{x,y\})  \\ & & 
+A_{x,y}\bar{B}_{x,y}\delta(C_1\cup\{x,y\})\delta(C_2)
+B_{x,y}\bar{B}_{x,y}\delta(C_1\cup C_2\cup\{x,y\}),
\end{eqnarray*}
and the two-component closure with $y\otr x\in C_1$ and $x\utr y\in C_2$
\[\scalebox{0.75}{\includegraphics{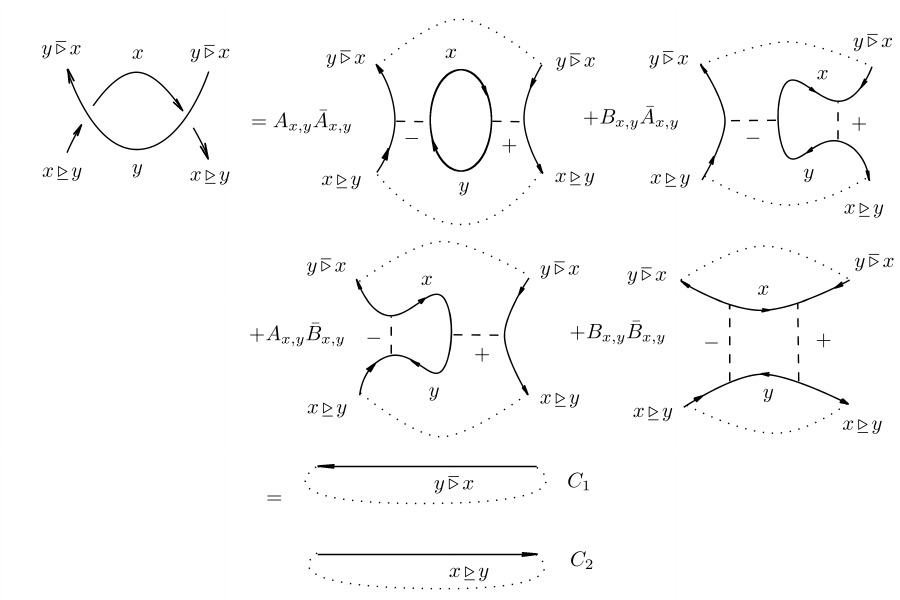}}\]
yields the requirement that 
\begin{eqnarray*}
\delta(C_1)\delta(C_2)
& = & A_{x,y}\bar{A}_{x,y}\delta(C_1\cup C_2)\delta(\{x,y\})
+B_{x,y}\bar{A}_{x,y}\delta(C_1\cup C_2\cup\{x,y\})\\ & & 
+A_{x,y}\bar{B}_{x,y}\delta(C_1\cup C_2\cup\{x,y\}) 
+B_{x,y}\bar{B}_{x,y}\delta(C_1\cup\{x\})\delta(C_2\cup\{y\}).
\end{eqnarray*}

Finally, the five equations in axioms (iii.i)-(iii.v) encode from the five ways of 
closing up an all-positive crossing Reidemeister III move portion of a Kauffman 
state. The trace diagrams on the two sides look like
\[\includegraphics{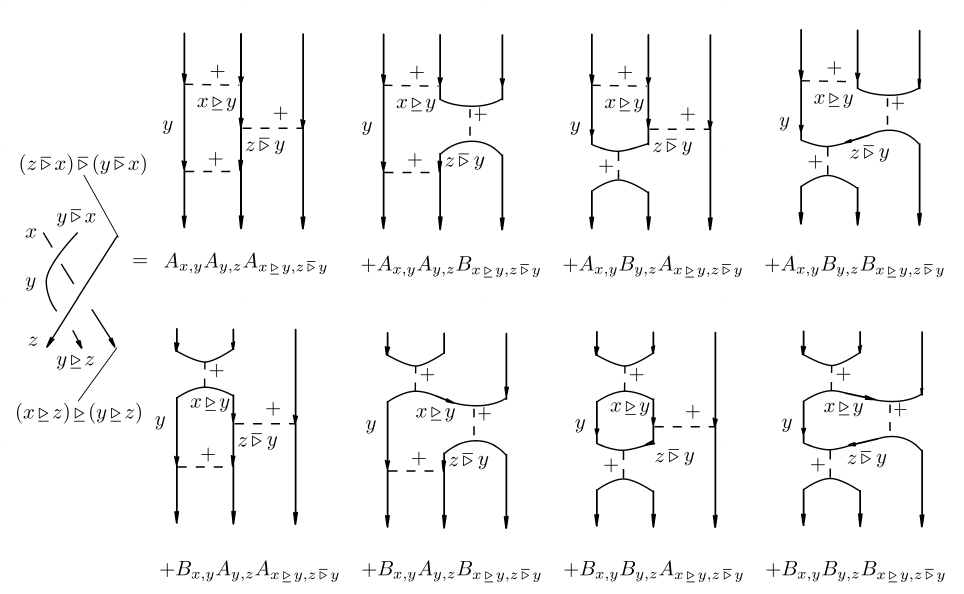}\]
and
\[\includegraphics{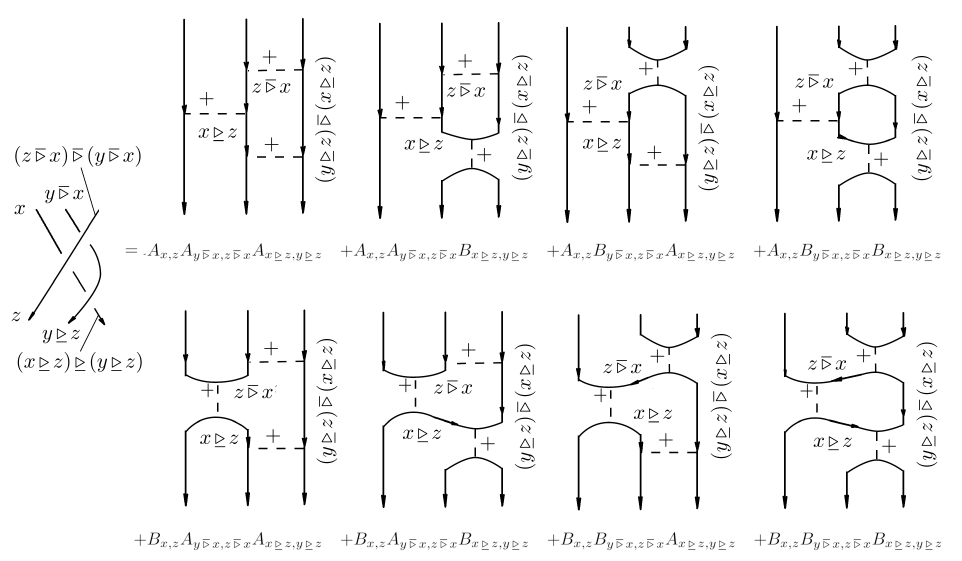}\]
and the possible ways of closing the components to obtain crossingless 
portions of Kauffman states are
\[\includegraphics{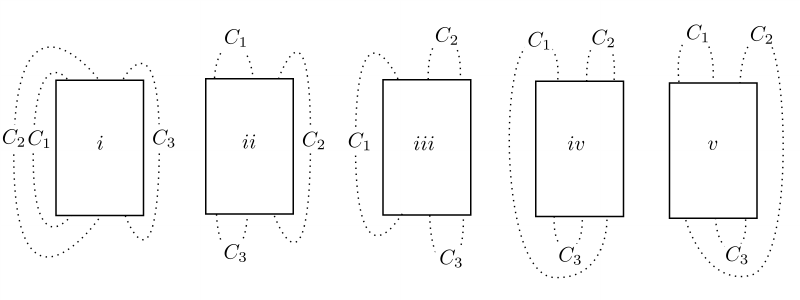}.\]

Hence, by construction we have our main result:

\begin{proposition}
Let $X$ be a finite biquandle, $R$ a commutative ring with identity,
($A,B,\bar{A},\bar{B}:X\times X\to R,w\in R^{\times},\delta:\mathcal{P}(X)\to R)$
a biquandle power bracket and $D$ an oriented link diagram with a choice of 
$X$-coloring and writhe $p-n$. Then the state-sum value $\beta_D$ obtained by 
summing over all Kauffman states of $D$ the products 
of smoothing coefficients times $\delta(C)$ for $C$ the set of biquandle colors
on each component of the state times the writhe correction factor $w^{n-p}$
is invariant under $X$-colored Reidemeister moves.
\end{proposition}

\begin{corollary}
The multiset of state-sum values over the set of $X$-colorings of an oriented 
knot or link is an invariant of knots and links.
\end{corollary}

As in the case of standard biquandle brackets, it is convenient to convert
the multiset into a polynomial format for ease of comparison. More precisely,
the multiset $\{\beta_1,\dots,\beta_n\}$ can be written as a formal polynomial
\[\sum_{k=1}^n u^{\beta_k}\]
in a formal variable $u$; then the multiplicity of each element becomes a
coefficient in $\mathbb{Z}$; then evaluation of the polynomial at $u=0$
yields the cardinality of the multiset.

\section{\large\textbf{Examples}}\label{E}

In this section we collect a few examples.

\begin{example}\label{ex1}
Biquandle brackets as defined in earlier papers such as \cite{NOR,NO}
are examples of biquandle power brackets with constant $\delta$ functions, 
specifically
\[\delta(C)=-A_{x,y}B_{x,y}^{-1}-A_{x,y}^{-1}B_{x,y},\quad C\in\mathcal{P}(X)\]
where this value is the same for all $x,y\in X$. However, for biquandle 
power brackets we no longer require the skein coefficients to be invertible 
elements of the coefficient ring, and hence we can have examples of 
biquandle power brackets with constant delta functions which were missing 
from the previous set of biquandle brackets.
For example, let $X=\{1,2\}$ be the biquandle defined by the operation tables
\[
\begin{array}{r|rr}\utr & 1 & 2 \\ \hline 1 & 2 & 2 \\ 2 & 1 & 1 \end{array}
\quad
\begin{array}{r|rr}\otr & 1 & 2 \\ \hline 1 & 2 & 2 \\ 2 & 1 & 1 \end{array}
\]
and let $R=\mathbb{Z}_4$. Then our \texttt{python} computations say that
the coefficient tables
\[
\begin{array}{r|rr}A & 1 & 2 \\ \hline 1 & 1 & 3 \\ 2 & 1 & 1 \end{array}\quad
\begin{array}{r|rr}B & 1 & 2 \\ \hline 1 & 0 & 0 \\ 2 & 2 & 0 \end{array}\quad
\begin{array}{r|rr}\bar{A} & 1 & 2 \\ \hline 1 & 0 & 3 \\ 2 & 0 & 2 \end{array}\quad
\begin{array}{r|rr}\bar{B} & 1 & 2 \\ \hline 1 & 3 & 0 \\ 2 & 3 & 1 \end{array}
\]
with $w=3$ and $\delta(S)=3$ for $S\in \{\emptyset, \{1\}, \{2\}, \{1,2\}\}$
defines a biquandle power bracket.
\end{example}

\begin{example}
With the same biquandle from Example \ref{ex1}, our computations say that
the coefficient tables over $\mathbb{Z}_5$
\[
\begin{array}{r|rr}A & 1 & 2 \\ \hline 1 & 1 & 4 \\ 2 & 0 & 1 \end{array}\quad
\begin{array}{r|rr}B & 1 & 2 \\ \hline 1 & 3 & 2 \\ 2 & 4 & 3 \end{array}\quad
\begin{array}{r|rr}\bar{A} & 1 & 2 \\ \hline 1 & 2 & 0 \\ 2 & 4 & 4 \end{array}\quad
\begin{array}{r|rr}\bar{B} & 1 & 2 \\ \hline 1 & 0 & 2 \\ 2 & 3 & 2 \end{array}
\]
with $w=2$ and $\delta(S)=4$ for $S\in \{\emptyset, \{1\}, \{2\}, \{1,2\}\}$
defines a biquandle power bracket.

We computed the value of the invariant for each of the prime links with up to 
7 crossings as listed at the Knot Atlas \cite{KA}; the results are collected 
in the table.
\[
\begin{array}{r|l}
\Phi_{X}^{\beta}(L) & L \\ \hline
2u+2u^2 & L6a2, L6a3 \\
2u+2u^3 & L2a1, L7a5, L7a6 \\
2u+2u^4 & L4a1, L6a1, L7a2, L7n1 \\
4u & L5a1, L7a1, L7a3, L7a4, L7n2 \\
2u+6u^4 & L7a7 \\ 
6u+2u^4 & L6a4, L6a5, L6n1
\end{array}
\]
\end{example}

Of course, we want biquandle power brackets with non-constant $\delta$
functions. 

\begin{example}
Let $X$ be the biquandle with operation tables
\[
\begin{array}{r|rrr}
\utr & 1 & 2 & 3 \\ \hline
1 & 2 & 2 & 1 \\
2 & 1 & 1 & 2 \\
3 & 3 & 3 & 3
\end{array}\quad
\begin{array}{r|rrr}
\otr & 1 & 2 & 3 \\ \hline
1 & 2 & 2 & 2 \\
2 & 1 & 1 & 1 \\
3 & 3 & 3 & 3
\end{array}.
\]
Then our \texttt{python} computations say that the coefficient tables over
$\mathbb{Z}_5$ 
\[
\begin{array}{r|rrr}
A & 1 & 2 & 3\\ \hline 
1 & 2 & 2 & 4 \\
2 & 4 & 3 & 0 \\
3 & 2 & 0 & 2
\end{array}\quad
\begin{array}{r|rrr}
B & 1 & 2 & 3\\ \hline 
1 & 0 & 4 & 0 \\
2 & 3 & 1 & 2 \\ 
3 & 2 & 3 & 3
\end{array}\quad
\begin{array}{r|rrr}
\bar{A} & 1 & 2 & 3\\ \hline 
1 & 4 & 0 & 0 \\
2 & 2 & 0 & 1 \\
3 & 4 & 1 & 3
\end{array}\quad
\begin{array}{r|rrr}
\bar{B} & 1 & 2 & 3\\ \hline 
1 & 1 & 3 & 4 \\
2 & 1 & 2 & 3 \\
3 & 0 & 3 & 2
\end{array}
\]
with $w=3$ and 
\[\delta(S)=\left\{ \begin{array}{ll}
4 & S\in \{\{1\}, \{2\}, \{1,2\}\}\\
0 & S\in \{\emptyset, \{3\},\{1,3\},\{2,3\},\{1,2,3\}\}
\end{array}\right.
\]
define a biquandle power bracket. The invariant values for prime
links with up to seven crossings are 
\[
\begin{array}{r|l}
\Phi_X^{\beta}(L) & L \\ \hline
1+2u+2u^2 & L2a1, L7a5, L7a6 \\ 
1+2u+2u^3 & L6a2, L6a3 \\
5+4u & L5a1, L7a1, L7a3, L7a4, L7n2 \\
5+2u+2u^4 & L4a1, L6a1, L7a2, L7n1 \\
7+6u+2u^4 & L6a5, L6n1, L7a7 \\
19+6u+2u^4 & L6a4
\end{array}.
\]
In particular, this small example shows that these invariants are
proper enhancements of the counting invariant. 
\end{example}

\begin{example}\label{last}
Let $X$ be the biquandle with operation tables
\[
\begin{array}{r|rrrr}
\utr & 1 & 2 & 3 & 4 \\ \hline
1 & 2 & 2 & 2 & 2 \\
2 & 1 & 1 & 1 & 1 \\
3 & 3 & 3 & 4 & 4 \\
4 & 4 & 4 & 3 & 3 \\
\end{array}
\quad
\begin{array}{r|rrrr}
\otr & 1 & 2 & 3 & 4 \\ \hline
1 & 2 & 2 & 1 & 1 \\
2 & 1 & 1 & 2 & 2 \\
3 & 4 & 4 & 4 & 4 \\
4 & 3 & 3 & 3 & 3 \\
\end{array}
\]
and let $R=\mathbb{Z}_6$. Then our \texttt{python} computations
say that the coefficient tables 
\[
\begin{array}{r|rrrr}
A & 1 & 2 & 3 & 4 \\ \hline 
1 & 4 & 4 & 2 & 5 \\
2 & 2 & 4 & 3 & 1 \\
3 & 0 & 3 & 2 & 0 \\
4 & 3 & 5 & 5 & 4
\end{array}\quad
\begin{array}{r|rrrr}
B & 1 & 2 & 3 & 4 \\ \hline 
1 & 3 & 3 & 1 & 3 \\
2 & 3 & 5 & 4 & 0 \\ 
3 & 2 & 2 & 0 & 2 \\
4 & 2 & 0 & 5 & 4
\end{array}\quad
\begin{array}{r|rrrr}
\bar{A} & 1 & 2 & 3 & 4 \\ \hline 
1 & 3 & 5 & 2 & 1 \\
2 & 4 & 5 & 2 & 4 \\
3 & 2 & 0 & 1 & 0 \\
4 & 5 & 2 & 3 & 5
\end{array}\quad
\begin{array}{r|rrrr}
\bar{B} & 1 & 2 & 3 & 4 \\ \hline 
1 & 4 & 4 & 4 & 0 \\
2 & 3 & 4 & 1 & 2 \\
3 & 2 & 3 & 1 & 2 \\
4 & 1 & 4 & 4 & 0
\end{array}
\]
with $w=5$ and 
\[\delta(S)=\left\{ \begin{array}{ll}
3 & S\in \{\{1\}, \{2\}, \{1,2\}\}\\
4 & S\in \{\{3\}, \{4\}, \{3,4\}\} \\
0 & \mathrm{otherwise}
\end{array}\right.
\]
Then the invariant has the following values on the set of prime 
classical links with up to seven crossings:

\[\begin{array}{r|l}
\Phi_X^{\beta}(L) & L \\ \hline
2u^2 + 4u^3+2u^4 & L2a1, L6a2, L7a5, L7a6 \\
8+4u^3+4u^4 & L4a1, L5a1, L61a, L7a1, L7a2, L7a3, L7a4, L7n1, L7n2 \\
8u^3+8u^4 & L6a5, L6n1, L7a7 \\
48+8u^3+8u^4 & L6a4
\end{array}\]

\end{example}

\begin{example}
Using the biquandle and biquandle bracket in Example \ref{last}, let us 
illustrate the computation of the invariant for the link $L4a1.$ The 
reader can verify that there are sixteen $X$-colorings of this link, 
including for instance
\[\includegraphics{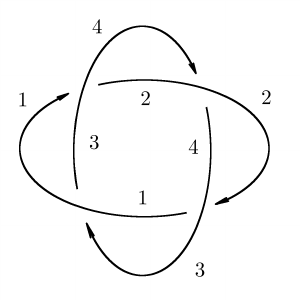}.\]
Selecting this coloring, let us compute the biquandle power bracket value. 
First, we find the set of $2^4=16$ Kauffman states, each of which has 
its associated coefficient list and $\delta$-value. In the interest 
of space, let us examine one of these states:
\[\includegraphics{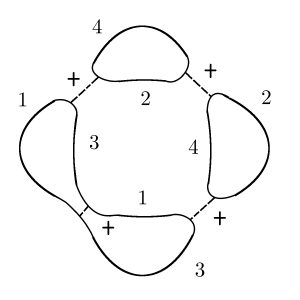}.\]
This state has coefficient list $A_{1,3}A_{4,2}A_{2,4}B_{3,1}$ and we have
$\delta$-values $\delta(\{1,3\})^2\delta(\{2,4\})$, so overall this state
contributes 
\[w^{-4}A_{1,3}A_{4,2}A_{2,4}B_{3,1}\delta(\{1,3\})^2\delta(\{2,4\})
=5^4(2)(5)(1)(2)(3)^2(4)=0.
\]
Multiplying the coefficients times delta values times writhe correction 
factor for each of the 16 Kauffman states and summing over the set, 
we obtain the biquandle power bracket value $0$ 
for this coloring. Repeating over the set of all colorings,
we obtain the multiset version of the invariant 
\[\{0,0,0,0,0,0,0,0,3,3,3,3,4,4,4,4\}\]
and converting to polynomial form yields $8+4u^3+4u^4$.
\end{example}

\section{\large\textbf{Questions}}\label{Q}

We end with a few questions for future research.

The biggest problem to solve going forward is how to find biquandle power 
brackets efficiently. For this project we have used a computer search
with our custom \texttt{python} code to find the few small examples
included. If $|X|=n$ and $|R|=m$, then the search space for all $X$-power 
bracket structures with $R$ coefficients has $4n^{2m+1}(2^n-1)^m$ elements,
so brute-force searches are inefficient. To find our examples, we first
found coefficients and $\delta$-values for the orbit sub-biquandles and
then filled in the other values with random values until solutions appeared.
As with the full brute-force search, this method's utility is limited to 
small-cardinality $X$ and $R$ cases. Moreover, our usual method of
using the axioms to propagate values though partial solutions cannot be used
in this case since the biquandle power bracket equations cannot be solved for
individual variables without the limiting assumption of invertibility of 
the variables.

In previous work such as \cite{FN}, biquandle bracket structures over infinite
rings such as $\mathbb{C}$ or $\mathbb{Z}[q^{\pm 1}]$ were found by ``lifting''
solutions over finite fields to the infinite case, e.g. going from 
$\mathbb{Z}_5$ to $\mathbb{C}$ by replacing $2$ with $i$ (since in 
$\mathbb{Z}_5$, $2^2=4=-1$) and verifying the biquandle bracket equations
for the resulting guess. Our attempts to repeat this for biquandle power 
brackets have not yielded any successes, suggesting that new methodology
is needed to find this type of solution as well.

Other questions of interest include generalizing biquandle power brackets
to other structures and cases. We note that unlike many classical knot 
invariants, biquandle power brackets do not extend to virtual knot invariants
by merely ignoring the virtual crossings in general; extending these 
invariants to the virtual case will take place in a sequel paper.

Understanding the structure of biquandle power brackets is of interest; in the
earlier case of standard biquandle brackets, it is known that such brackets
can be cohomologous to each other via biquandle 2-coboundaries. What is the
algebraic structure, if any, of the set $\mathcal{BPB}(X,R)$ of biquandle 
power brackets over a given biquandle $X$ and coefficient ring $R$, and how 
does it relate to that of the standard case? What is the relationship between 
$\mathcal{BPB}(X,R)$ and $\mathcal{BPB}(X',R)$ or $\mathcal{BPB}(X,R')$ for 
sub-biquandles $X'\subset X$  or extension rings $R'\supset R$?

\bibliography{ng-sn2}{}
\bibliographystyle{abbrv}

\bigskip

\noindent\parbox{3in}{
\textsc{Department of Mathematics\\ 
Izmir Institute of Technology\\
Urla 35433 \\
Izmir, Turkey}}
\quad
\parbox{3in}{
\textsc{Department of Mathematical Sciences \\
Claremont McKenna College \\
850 Columbia Ave. \\
Claremont, CA 91711}}

\end{document}